\documentclass[a4paper]{article}

\usepackage{graphicx}

\usepackage{mathrsfs}
\usepackage{amsmath}			
\usepackage{amssymb}		
\usepackage{csquotes}
\usepackage{hyperref}
\hypersetup{colorlinks=true, linkcolor=black, citecolor=black}
\usepackage{dsfont}
\usepackage{mathtools}
\usepackage{algorithm}
\usepackage{algorithmic}
\usepackage{float}

\newcommand\norm[1]{\lVert#1\rVert}
\newcommand{\R}{\mathbb{R}}
\newcommand{\mc}[1]{\mathcal#1}

\newcommand{\trace}{\mathop\mathrm{tr}}

\title{ Adaptive Gaussian Process Regression for Bayesian inverse problems\thanks{This work has been supported by Bundesministerium für Bildung und Forschung -- BMBF, project number 05M20ZAA (siMLopt) and by the Deutsche Forschungsgemeinschaft (DFG, German Research Foundation) -- 436400679.
}}

\author{Paolo Villani\thanks{Zuse Institut Berlin, \{weiser,villani\}@zib.de} \and
        Jörg Unger\thanks{Bundesanstalt für Materialforschung und -prüfung, joerg.unger@bam.de}
        \and Martin Weiser\footnotemark[2]}

\begin{document}

\maketitle

\begin{abstract}
We introduce a novel adaptive Gaussian Process Regression (GPR) methodology for efficient construction of surrogate models for Bayesian inverse problems with expensive forward model evaluations. An adaptive design strategy focuses on optimizing both the positioning and simulation accuracy of training data in order to reduce the computational cost of simulating training data without compromising the fidelity of the posterior distributions of parameters. The method interleaves a goal-oriented active learning algorithm selecting evaluation points and tolerances based on the expected impact on the Kullback-Leibler divergence of surrogated and true posterior with a Markov Chain Monte Carlo sampling of the posterior. The performance benefit of the adaptive approach is demonstrated for two simple test problems.
\end{abstract}



\pagestyle{myheadings}
\thispagestyle{plain}
\markboth{M. WEISER AND P. VILLANI}{GAUSSIAN PROCESSES FOR INVERSE PROBLEMS}

\section{Introduction}

The inverse problem of inferring the posterior probability of parameters $p\in\R^d$ in a forward model $y(p)$ from measurements $y^m\in\R^m$ is often addressed by sampling with Markov Chain Monte Carlo (MCMC) methods~\cite{GreenLatuszynskiPereyraRobert2015}. The large number of forward evaluations required for a faithful representation of the posterior density renders this inapplicable in case of computationally expensive forward models such as large finite element (FE) simulations. The forward model is thus often replaced by a fast surrogate model when sampling the posterior. Here, we focus on the efficient construction of Gaussian Process Regression (GPR) surrogates.

Surrogate models are learned from values $y(p_i)$ at specific evaluation points $p_i$ as training data. The accuracy of the resulting surrogate depends on the number and position of the sample points. Constructing an accurate surrogate model can become computationally expensive when a large number of evaluations is required. Consequently, strategies for selecting near-optimal evaluation points have been proposed for various settings~\cite{RasmussenWilliams2006}. A priori point sets~\cite{Giunta,Queipo} are effectively supplemented by adaptive designs~\cite{Crombecq,Joseph,Lehmensiek,Sugiyama} selecting the most beneficial evaluation points $p_i$. 

When using FE simulations for computing training data, the evaluations of $y(p_i)$ are affected by discretization and truncation errors. The trade-off between accuracy and cost has been investigated using different low and high fidelity models~\cite{Nitzler}, and by an adaptive choice of  evaluation tolerances~\cite{SagnolHegeWeiser2016,SemlerWeiser2023,SemlerWeiser2024} in different settings.
Here, we extend~\cite{SemlerWeiser2023} from an offline training for maximum posterior point estimates to an interleaved posterior sampling and surrogate training driven by a goal-oriented approach.

\section{Gaussian Process regression}
Gaussian process regression is a regression technique which allows to approximate any function, naturally fits the Bayesian framework, and provides an uncertainty estimate of its prediction.

We consider a forward model $y : \R^d  \rightarrow \R^m $, which we assume to be a realisation of a Gaussian process $\mc G$ with mean $\mu_0 :  \R^d  \rightarrow \R^m$ and covariance kernel $k : \R^d \times \R^d  \rightarrow \R^{m\times m}$ to be defined later. 

For training points $(p_i,  y_i)_{i=1, \dots, s }$ with $y_i\approx y(p_i)$ of accuracy $\tau_i\ge 0$, we are interested in a prediction of $y_{s+1} \approx y(p_{s+1})$ for any $p_{s+1}$. The GPR posterior covariance block matrix is $\Gamma = (K^{-1} + T^{-2})^{-1}\in \R^{m(s+1) \times m(s+1)}$ with prior covariance blocks $K_{ij} = k(p_i,p_j)$ and formally likelihood covariance $T=\mathrm{diag}(\tau_1 I,\dots,\tau_s I,\infty I)$. The GPR posterior mean is $\bar Y = \Gamma(K^{-1}M_0 + T^{-2}Y)$ with $Y=(y_1,\dots, y_s, 0)$. Then, the GPR prediction is the marginal normal distribution $y_{s+1}\sim\mathcal{N}(\bar Y_{s+1},\Gamma_{s+1,s+1})$. As $p_{s+1}\in\Omega$ is arbitrary, this defines mean $\bar y:\Omega\to\R^m$ and covariance $\Gamma:\Omega\to\R^{m\times m}$ on the whole parameter space. We refer to~\cite{RasmussenWilliams2006,SemlerWeiser2023} for a more detailed exposition.

\section{Bayesian surrogate-based parameter identification}
We consider the forward model $y : \Omega \subset \R^d  \to \R^m$, which cannot be evaluated directly, but can be approximated through a numerical procedure $y_\tau$ with arbitrary precision in exchange of computational work: We assume that for any $\tau>0$, we obtain an evaluation $y_\tau(p)\sim\mc N(y(p),\tau I)$, with cost $W_\tau$.

We assume measurements $y^m$ to be random variables generated by a linear additive Gaussian noise model
\begin{equation} \label{eq:meas-mod}
    y^m = y(p) + \eta
\end{equation}
with $\eta \sim \mc N (0, \Sigma_l)$. 
For simplicity, we consider a diagonal covariance structure $\Sigma = \text{diag}( \sigma_1, \dots \sigma_n)$, corresponding to independent noise components. 
The conditional distribution of the measurements is then $y^m \mid p \sim \mathcal N (y(p), \Sigma_l )$,
\[
\pi (y^m \mid p) = (2\pi)^{-m/2} \det(\Sigma_l)^{-1/2} \exp \Big (-\frac{1}{2}\norm{y^m - y(p)}_{\Sigma_l^{-1}}^2 \Big )
\] 
is the likelihood of the problem. Evaluating the likelihood requires evaluating the forward model $y$, which we assume to be computationally expensive.

To reduce costs, we assume that $y$ is a realisation of a GP, and introduce a GP surrogate model $\mc G$ of predictive mean $\bar y$ and variance $\Gamma$. For simplicity, we consider a surrogate with independent output components, i.e. diagonal covariance $\Gamma(p)$. 
The training points for this GP are given by numerical evaluations $y_{\tau_i}(p_i)$ of the forward model. These points and the corresponding evaluation tolerances $\tau_i$ form the training design $\mc D$. We postpone the question of how to build training designs to the next section.

\begin{figure}[ht] 

\begin{centering}
\includegraphics[width = 360pt]{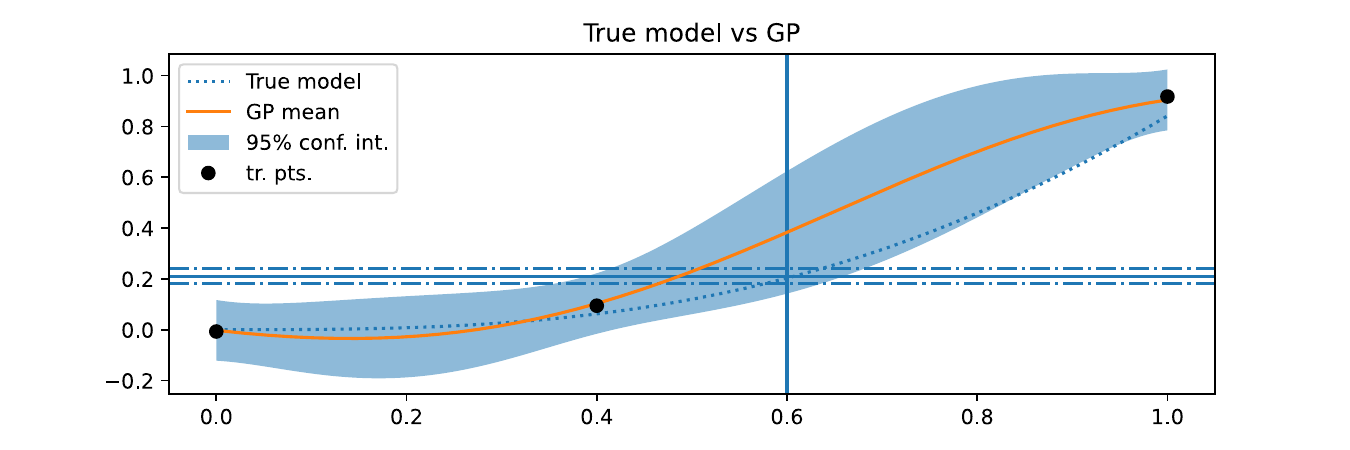}
\includegraphics[width = 360pt]{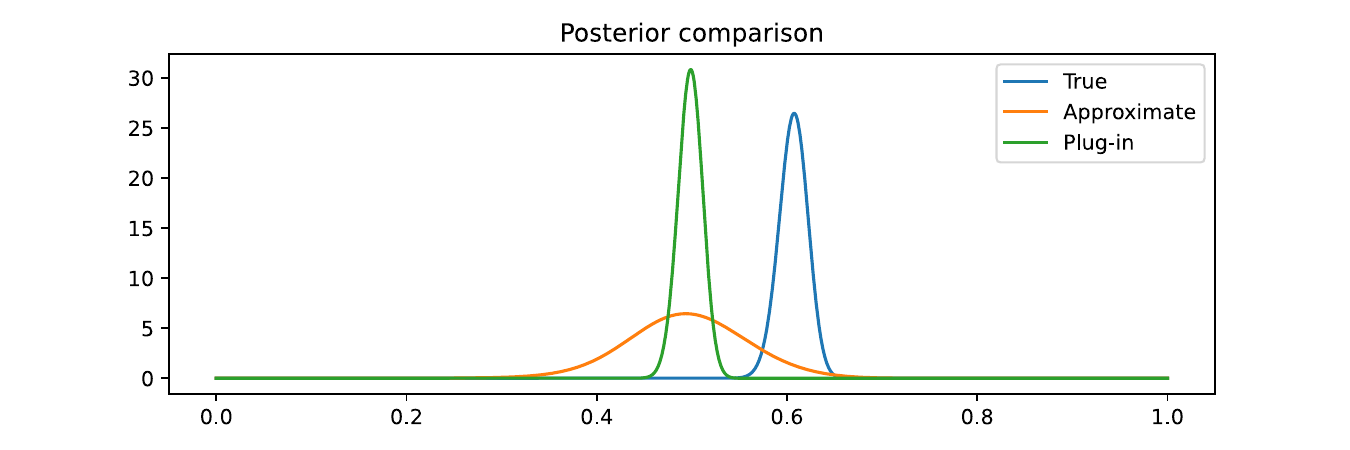}
\end{centering}

\caption{Impact of the likelihoods~\eqref{eq:plug-in-likelihood} and~\eqref{eq:full-likelihood} on the posterior for an illustrative inverse problem problem with forward model $y(p) = p^2 sin(p)$ and uniform prior on parameter space $[0,1]$. The horizontal lines show the actual measurement and the $2\sigma$ range of measurement noise. The marginal likelihood~\eqref{eq:full-likelihood} is wider due to including the GP variance, and avoids overconfident posteriors.} 
\label{fig:likelihoods}
\end{figure}

To evaluate the likelihood, we could substitute the forward model with the mean estimate $\bar y$, obtaining 
\begin{equation}\label{eq:plug-in-likelihood}
\pi_{\text{plug-in}}(y^m\mid p, \bar y) = (2\pi)^{-m/2} \det(\Sigma_l)^{-1/2}\exp \Big (-\frac{1}{2}\norm{y^m - \bar y(p)}_{\Sigma_l^{-1}}^2 \Big ).
\end{equation}
This, from a decision-theoretic point of view, corresponds to the minimisation of the $L^1$ loss \cite{JarvenpaaGutmannVehtariMarttine2021}, but ignores the uncertainty estimate given by the predictive variance: since $y$ is assumed to be a realisation of $\mc G$, the measurement noise model \eqref{eq:meas-mod} becomes $y^m = \mc G(p) + \eta$.
Marginalizing over GP realizations results in a different conditional distribution of the measurements
$y^m \mid p, \mc G \sim \mathcal N (\bar y(p), \Sigma_l + \Gamma (p) )$
and in a marginal likelihood: 
\begin{equation}\label{eq:full-likelihood}
\pi_\mc D(y^m\mid p, \mc D) = (2\pi)^{-m/2}\det\left ( \Sigma_l +\Gamma(p) \right )^{-\frac{1}{2}} \exp \left ( -\frac{1}{2}\norm{y_m - \bar y (p)}_{\left (\Sigma_l + \Gamma(p) \right ) ^{-1}}^2 \right ),
\end{equation}
see, e.g.,~\cite{BaiTeckentrupZygalakis2023}.
Note that the conditional distribution is still Gaussian due to the normality of both the noise and the GP. Moreover, the likelihood $\pi_\mc D$ is closely related to the $L^2$ loss~\cite{JarvenpaaGutmannVehtariMarttine2021, SinsbeckNowak}. Including the GP variance into the likelihood can be important for avoiding overconfident yet wrong posterior approximations by surrogated forward models, see Fig.~\ref{fig:likelihoods} for an illustration.

By adopting a Bayesian point of view, we express prior belief on the parameter by assigning a prior distribution $\pi(p)$. Then, by Bayes' theorem, we obtain a true posterior distribution
\begin{equation} \label{eq:true_post}
    \pi(p \mid y^m ) =
    \frac{\pi(p) \ \pi(y^m \mid p )}{\pi(y^m)} ,
\end{equation}
corresponding to the true likelihood $\pi(y^m \mid p)$ and an approximate posterior 
\begin{equation} \label{eq:approx_post}
    \pi(p \mid y^m, \mc D) =
    \frac{\pi(p) \ \pi(y^m \mid p, \mc D )}{\pi(y^m \mid \mc D)} ,
\end{equation}
corresponding to the likelihood approximation $\pi(y^m \mid p,\mc D)$ as given in~\eqref{eq:plug-in-likelihood} and~\eqref{eq:full-likelihood}, respectively.

In both cases, the normalising constant $\pi(y^m)$ or $\pi(y^m\mid \mc D)$, respectively, will not be computationally available, as it requires integration over the parameter space $\Omega$: fortunately, it is not needed for posterior sampling by Markov-Chain Monte Carlo (MCMC) methods.

\section{Posterior-oriented surrogate model}

As in \cite{SinsbeckNowak}, we do not aim  at building a surrogate which is globally accurate on the whole parameter space $\Omega$, but at finding a design $\mc D$ such that the approximate posterior is accurate, i.e. $\pi( p \mid y^m) \approx \pi( p \mid y^m, \mc D)$. Repeatedly selecting training points randomly sampled from $\pi(p \mid y^m, \mc D)$, updating $\mc G$ and then iterating is sufficient for convergence of $\pi(p \mid y^m, \mc D)$ to $\pi(p \mid y^m ) $ in the Hellinger metric~\cite{BaiTeckentrupZygalakis2023}. Here, we also aim at finding a design $\mc D$ which incurs a small computational cost of evaluating training data $y_{\tau_i}$.

We measure the deviation of the surrogated and the true posterior densities by the Kullback-Leibler (KL) divergence
\begin{align}
    D_{\rm KL}\left (\pi(\cdot \mid y_m) \mid \pi(\cdot \mid y_m, \mc D)  \right ) &= \mathbb E_{\pi(p \mid y^m)} \left [ \log  \frac{\pi(p \mid y_m)}{\pi (p \mid y_m, \mc D)}  \right ]\notag \\
    &= \int_{\Omega} \pi(p \mid y_m) \log \frac{\pi(p \mid y_m)}{\pi(p \mid y_m,\mc D)} \ dp. \label{eq:error-quantity}
\end{align}
Since computing the KL divergence requires evaluating the full model, we derive a numerical approximation which relies on the surrogate only. Using the marginal likelihood~\eqref{eq:full-likelihood} and the posteriors~\eqref{eq:true_post} and~\eqref{eq:approx_post}, their logarithmic ratio can be written as
\[
\log  \frac{\pi(p \mid y^m)}{\pi (p \mid y^m, \mc D)}  = 
\log  \frac{\pi( y^m \mid p)}{\pi_\mc D ( y^m \mid p,  \mc D)}  -
\log  \frac{\pi(y^m)}{\pi ( y^m \mid \mc D)} .
\]
The first term, the logarithmic ratio of true and surrogated likelihood, equals
\begin{align*}  
\log&  \frac{\pi( y^m \mid p)}{\pi_\mc D ( y_m \mid p,  \mc D)}  \\
&= \frac{1}{2} \Bigg ( \log    \frac{\text{det}\left ( \Sigma_l + \Gamma(p)  \right )}{\text{det} \left (\Sigma_l \right )} 
- \norm{ y (p) - y^m}_{\Sigma_l^{-1}} ^2 +  \norm{ \bar y (p) - y^m}_{\left (\Sigma_l + \Gamma(p)\right )^{-1}}^2 \Bigg ).
\end{align*}
As $\Sigma_l^{-1} - \left (\Sigma_l + \Gamma (p) \right )^{-1} \succeq 0$, we can upper bound the difference between norms by
\begin{align*}   
    -\norm{ y (p) &- y^m}_{\Sigma_l^{-1}} ^2 +  \norm{ \bar y (p) - y^m}_{(\Sigma_l + \Gamma(p))^{-1}}^2  \\
    &\leq -\norm{ y (p) - y^m}_{\Sigma_l^{-1} } ^2 +  \norm{ \bar y (p) - y^m}_{\Sigma_l^{-1}}^2  \\
    &= -\norm{ y (p) - \bar y(p)}_{\Sigma_l^{-1} } ^2 - 2  \left ( \bar y (p) - y^m \right ) ^T \Sigma_l^{-1} \left ( y (p) - \bar y (p) \right ) .
\end{align*}

By assuming that $y$ is a realisation of $\mc G$, $\mathbb E \left [ \left ( y^{(i) }(p) - \bar y^{(i) }(p) \right ) ^2 \right ] = \Gamma^{(i,i) }(p)$  and therefore $\norm{ y (p) - \bar y(p)}_{\Sigma_l^{-1} } ^2 \approx \mathop\mathrm{tr} \left(\Sigma_l^{-1}\Gamma(p) \right )$ hold. Defining $v = \Sigma_l^{-1}  \sqrt{ \mathop\mathrm{diag} \left(\Gamma(p) \right )} \in \R^m$, we obtain
\begin{align*}
    -\norm{ y (p) - y^m}_{\Sigma_l^{-1}} ^2 +  \norm{ \bar y (p) - y^m}_{(\Sigma_l + \Gamma(p))^{-1}}^2  
    \lesssim - \mathop\mathrm{tr} \left(\Sigma_l^{-1}\Gamma(p) \right )+2 \left | \bar y (p) - y^m  \right | ^T v .
\end{align*}
We therefore define the local error quantity
\begin{align} \label{eq:num-logdiff}
e_\mc D (p) & := \frac{1}{2} \left ( \log\det(I+\Sigma_l^{-1}\Gamma(p)) - \trace\left(\Sigma_l^{-1}\Gamma(p)(I-\Sigma_l^{-1}) \right ) 
+2 \left | \bar y (p) - y^m \right | ^T v \right ) \\
& \gtrsim \log  \frac{\pi( y^m \mid p)}{\pi_\mc D ( y_m \mid p,  \mc D)} \notag
\end{align}
as an approximate upper bound on the log ratio of true and surrogated likelihood.

By optimistically assuming that the normalisation factors are similar independent of the design $\mc D$, and thus $\log  \frac{\pi(y^m)}{\pi ( y^m \mid \mc D)}  \approx 0 $, we substitute~\eqref{eq:num-logdiff} into~\eqref{eq:error-quantity} and obtain the global error quantity
\begin{align} \label{eq:target}
    E (\mc D )= \int_\Omega e_\mc D (p) \pi(p \mid y^m) \ d p.
\end{align}

To create an optimal surrogate model, we aim at a training design $\mc D$ which minimises $E(\mc D)$ under a computational work constraint. By denoting the computational work needed to realize $\mc D$ by  $W(\mc D)$ , for a given budget $W$ we aim at solving the optimisation problem
\begin{equation} \label{eq:doe}
    \min_{\mc D}  E(\mc D) \ \text{ subject to } \ W(\mc D) \leq W.
\end{equation}

\section{Sequential design of experiments}
It is far from trivial to predict  a priori how design choices impact the error quantity $E$, especially when a large budget $W$ is available or the initial surrogate is unreliable. Fortunately, an exact solution of~\eqref{eq:doe} is not needed -- an approximate solution will do, even if it yields a slightly less efficient design. We follow~\cite{SemlerWeiser2023, SemlerWeiser2024} and adopt a greedy sequential approach, where the budget $W=\sum_{j=1}^J \Delta W_j$ is partitioned and sequentially spent. 

We start from an initial design $\mc D_0$ and then, for $j=1,\dots ,J$, aim at solving
\begin{equation} \label{eq:incremental-doe}
    \min_{\mc D_{j} \le \mc D_{j-1}} E(\mc D_{j}) \quad \text{s.t.} 
    \quad W(\mc D_{j} \mid \mc D_{j-1}) \le \Delta W_j.
\end{equation}
We write $\mc D \le \mc D_{j-1}$ for any design $\mc D$ which refines $\mc D_{j-1}$ in the sense that it includes all evaluation points $p_i$ contained in $\mc D_{j-1}$ with lesser or equal tolerances $\tau_i$. We write $W(\mc D \mid \mc D_{j-1}) = W(\mc D) - W(\mc D_{j-1})$ for the work needed to obtain $\mc D$ from $\mc D_{j-1}$.

Even this sequential formulation is highly non-linear and non-convex. An accurate solution would require a considerable amount of computational work, possibly exceeding the savings in computational budget possible with a better design. Consequently, we adopt the heuristic approach of separating the selection of new candidate evaluation points from the optimisation of the evaluation tolerances. In the latter, we also decide about the actual inclusion of the new points in the training set.

\emph{Candidate points.} We choose points where spending computational budget is likely to reduce the error most. In order to do so, we look at the sensitivity of the global error $E$ with respect to a reduction of training error at a candiate position $p'$~\cite{SemlerWeiser2024}. This is given by
\begin{align} 
    \frac{ d E ({\mc D}) }{d W(p')} 
    &= 
    \int_\Omega \frac{d e_{{\mc D}} (p) }{d W(p')} \pi_{\mc D}(p\mid y_m) \, dp \notag \\
    &=  \int_\Omega \frac{d e_\mc D(p)}{d \Gamma(p)}  \frac{d \Gamma(p)}{d \tau(p')}\bigg|_{\tau=\tau'}  \frac{d \tau(p')}{d W(p')}\bigg|_{\tau=\tau'}\pi_{\mc D}(p\mid y_m) \ dp , \label{eq:utility}
\end{align}
where the linearization tolerance $\tau '$ is the current GP standard deviation at point $p'$. We adopt \eqref{eq:utility} as a utility function and select local minimizers of  $\frac{ d E (\mc D_{j-1}) }{d W} $ as next candidate points.

The optimisation problem is solved approximately via a multistart pattern search. Quadrature is performed by Monte Carlo integration on samples $\mc S_j$ to be defined in Sec.~\ref{sec:sampling} below. This results in the numerical utility function
\[
\frac{ d E ({\mc D_{j-1}}) }{d W(p')} \approx \frac{1}{|\mc S_j|} \sum_{p \in \mc S_j} \frac{d e_{\mc D_{j-1}}(p) }{d W(p')} .
\]
If more than $c_j$ local maxima are found, the best $c_j$ ones are selected as candidates; if less are found, all of them are included. A larger number of candidates allows more points to be considered, but results in a harder accuracy optimisation problem.

\emph{Evaluation tolerances.} Let $\mc D_{j} = \big \{ (p_i^j,\tau_i^j) \mid i = 1, \dots, s_j \big \}$   be the set of training points at step $j$. By the selection of candidate points, $s_j \ge s_{j-1}$ and $p_i^j = p_i^{j-1}$ for $i=1,\dots, s_{j-1}$ hold.

Optimal tolerances $\tau_i^j$ are given by the solution of~\eqref{eq:incremental-doe} as a function of the tolerances. In order to be able to solve the problem, we ignore the shifts in the mean $\bar y$ as they cannot be predicted before evaluating the model. Consequently, we only consider the impact of evaluation tolerances on the predictive variance and, for evaluation tolerances $\tau^j = (\tau_1^j, \dots, \tau_{s_j}^j)$, write $E( \tau^j)$. As already spent computational budget cannot be recovered by forgetting previously acquired information, we impose the constraint $\tau_i^j \le \tau_i^{j-1}$ for $i=1,\dots,s_{j-1}$.

This results in the problem
\begin{equation} \label{eq:accuracy-doe}
     \min_{\tau^j \in \mc T_j} E(\tau^j) \quad \text{subject to} \quad W_{\tau^j | \mc D_{j-1}}\leq \Delta W_j, 
\end{equation}
where the set of admissible tolerances is
\[
\mc T_j = \{ (\tau _1, \dots, \tau _{s_j}) \in ( \R ^+ \cup \{ +\infty \} ) ^{s_j} \mid \tau_i \leq \tau_i^{j-1} \text{ for } i \leq s_{j-1}  \}.
\]
If after optimization $\tau_i^j = +\infty$ holds for some $i>s_{j-1}$, $p_i^j$ is excluded from the training set.

Before we can numerically solve the problem, we need to notice that computational costs are not available before the evaluation is performed, such that we need to resort to a priori work models. Following~\cite{SemlerWeiser2023,WeiserGhosh2018}, we make use of established a priori asymptotic estimates for finite elements of degree $r$ in space dimension $l$ and an optimal solver such as multigrid, and define
\begin{equation} \label{eq:work-model}
    W(\tau) = \tau^{-l/r} .
\end{equation}
This estimate is asymptotic for $\tau \rightarrow 0$.  Consequently, despite being inaccurate for low-accuracy evaluations, it is usually accurate for the expensive high-accuracy ones.

Problem \eqref{eq:accuracy-doe} is solved by multistart gradient descent with projection and backtracking linesearch. The integral in $E$ is approximated again by Monte Carlo integration on the samples $\mc S_j$, resulting in a numerical objective 
\[
E(\tau^j) \approx \frac{1}{|\mc S_j|} \sum_{p \in \mc S_j} e_{\tau^j}(p).
\] 
To implement gradient descent with projection, we adopt the coordinate change 
\[
\tau^j = \left (\tau _1, \dots, \tau _{s_j}\right ) \mapsto \left ( \tau _1^{-l/r}, \dots, \tau _{s_j}^{-l/r} \right ) = W^j,
\]
such that the constraint in~\eqref{eq:accuracy-doe} becomes linear, transforming the set of admissible tolerances $\mc T^j$ into a simplex and enabling efficient projection.

\section{Solution of the inverse problem} \label{sec:sampling}

The previous sections established the inverse problem~\eqref{eq:true_post} and the sequential approach \eqref{eq:incremental-doe} to surrogate model training. Similar to~\cite{WangBroccardo2020}, we combine them to an interleaved strategy given as pseudocode in Alg.~\ref{alg:interleaved}.

Both the global error quantity~\eqref{eq:target} and the utility function~\eqref{eq:utility} require integration with respect to the posterior $\pi(p \mid y^m)$. We perform the integration through an MCMC sampling of the posterior, which is is at the same time the ultimate goal of the inversion. 

We start with an empty sample chain $\mc S_0 = \emptyset$. At iteration $j$, we draw a number $n_j$ of samples form $\pi(p\mid y^m, \mc D_{j-1}) $, append them to $\mc S_{j-1}$, and remove the oldest $h_j < n_j$ elements of the chain, as they have been drawn from a less accurate posterior approximation. This results in the sample chain $\mc S_j$, which is used to evaluate the integrals involved in the training problem \eqref{eq:incremental-doe} at step $j$.

As the sample size $|\mc S_j|$ may be is too large for an efficient evaluation of the integrals in~\eqref{eq:target} and~\eqref{eq:utility}, we use a sufficiently large randomly extracted subset of $\mc S_j$ instead of the whole chain for Monte Carlo integration.

When the  computational budget is exhausted, the training of the surrogate model terminates. A last round of samples is added to the chain, obtaining the final set of samples from the posterior.
\begin{algorithm}[]
\caption{Surrogate-based Bayesian inversion} \label{alg:interleaved}
\begin{algorithmic}[1]
\REQUIRE $\mc D_0$ initial design, $W$ budget
\STATE $\mc S_0 \leftarrow \emptyset$
\STATE $W _\mc D \leftarrow 0$
\STATE $j \leftarrow 1$
\WHILE{ $W_\mc D \leq W $ }
\STATE \textbf{decide:} $n_j$ samples to draw, $h_j$ samples to remove
\STATE remove $h_j$ samples from $\mc S_{j-1}$
\STATE draw $n_j$ samples $\mc S$ from $\pi(p \mid y^m, \mc D)$
\STATE $\mc S_j \leftarrow \mc S_{j-1} \cup \mc S$
\STATE \textbf{decide:} $\Delta W_j$ iteration budget, $c_j$ number of candidates
\STATE obtain $c_j$ candidates
\STATE optimize accuracies $\tau^j$,  update $\mc D$
\STATE evaluate forward model for decreased tolerances
\STATE $W_\mc D\leftarrow W_\mc D + \Delta W_j $
\STATE $j \leftarrow j + 1$
\ENDWHILE
\STATE draw $n_j$ samples $\mc S$ from $\pi(p \mid y^m, \mc D)$
\STATE $\mc S_j \leftarrow \mc S_{j-1} \cup \mc S$
\end{algorithmic}
\end{algorithm}

\section{Numerical experiments}
We present two illustrative experi\-ments based on a Python implementation of Alg.~\ref{alg:interleaved}, where GPR is implemented with PyTorch. We adopt a separable kernel with diagonal output structure and a Gaussian kernel as base~\cite{AlvarezRosascoLawrence2012}. The hyperparameters are tuned by marginal likelihood maximisation using PyTorch's Adam optimiser, with the kernel's correlation length scale constrained to $[0, 0.15]$.

As a benchmark, the results are compared with a non-adaptive space filling approach, Latin Hypercube Sampling, and the position-adaptive-only training strategy given by candidate point selection according to~\eqref{eq:utility}, i.e. all candidates are accepted and evaluated with a fixed accuracy.
For comparing the approaches, the approximation errors~\eqref{eq:error-quantity} are computed numerically with MCMC sampling utilising the true forward model. The implementation used for these examples is available at Zenodo\footnote{\url{https://zenodo.org/doi/10.5281/zenodo.11066159}}.

\subsection{1D analytical experiment}
The first experiment is performed on a one-dimensional parameter space, with $m=2$ measurements. We consider an analytical forward model $y : \mathopen]0, 1 \mathclose[ \to \R^2$ given by
\begin{gather*}
    y(p) =  \left [
    \frac{1}{2} p + \frac{1}{2} p^2 \text{ exp}\left ( \frac{1}{3} \text{sin} ( 12 p - i) \right ) \right ]_{i = 0,1}.
\end{gather*}
This mimics the evaluation of a FE model on a 2D domain with quadratic elements, i.e. $l/r=1$. The discretization error is simulated via a zero mean Gaussian noise and the measurement likelihood is $\Sigma_l = 10^{-4} \mathrm{diag}( \frac{16}{9}, \frac{4}{9})$.

A budget of 500 is considered: at each iteration two candidate points are considered and a budget of 20 is assigned to each point. With the work model~\eqref{eq:work-model}, this results in a default tolerance of $0.05$ per point in the non-adaptive strategies and a total of 12 iterations.

The number $n_j$ of new samples added into $\mc S_j$ is gradually increased from 200 samples at the first iteration to 2000 in the last, according to $n_j = 200 + 1800 \left ( \frac{j}{12} \right )^2$. Similarly, the number of discarded samples ranges from 200 to 1000, with $h_1 = 0$ as in the first iteration the chain is empty, and $h_j = 200 + 800 \left ( \frac{j}{12} \right )^2$ for $j > 1$.

\begin{figure}
\begin{centering}
\includegraphics[width = \textwidth]{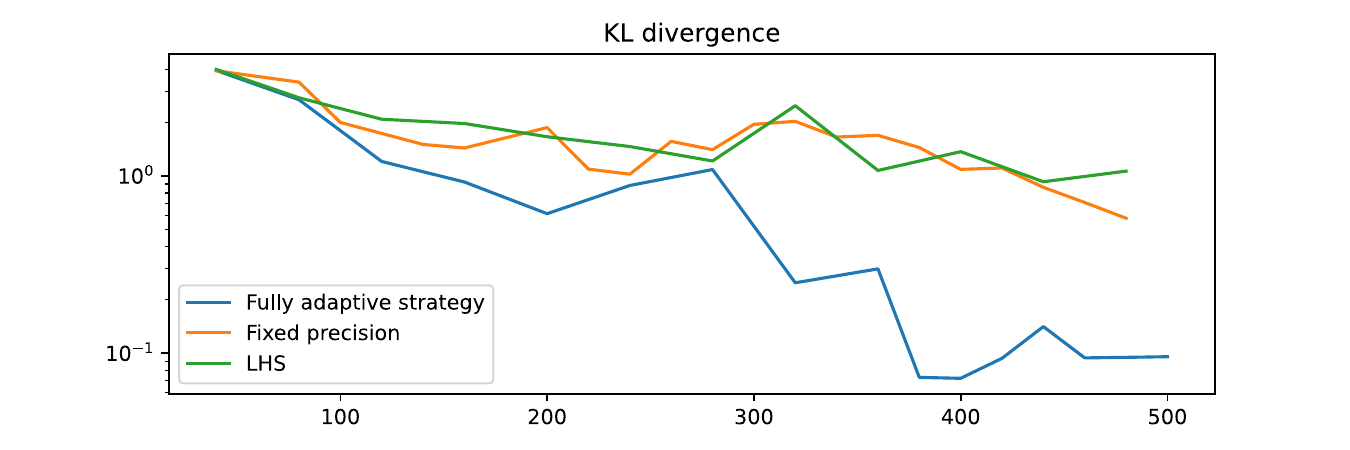}
\end{centering}
\caption{Kullback-Leibler divergence of surrogated posterior and true posterior for different training designs over the computational work spent in the 1D example.} 
\label{fig:1D}
\end{figure}

The obtained accuracies in terms of the Kullback-Leibler divergence between true posterior $\pi(p\mid y^m)$ and surrogated posterior $\pi(p\mid y^m,\mc D)$ are shown in Fig.~\ref{fig:1D}. Optimizing evaluation tolerances provides a significant performance improvement over both other strategies.

\subsection{2D analytical experiment}
The second experiment considers a parameter space of two dimensions and $m=3$ measurements. The forward model $y : \mathopen]-0.5, 0.5 \mathclose[^2 \to \mathbb R^3$ is again analytical, given by
\begin{align*}
    y(p) = 
     \bigg [ &\text{sin} ( 10k ) ( p_1 -  p_ 2 )  \text{ exp}\left ( \frac{1}{3} \text{sin} ( 8 p_2) \right ) \\
     &+ \text{cos} ( 10k ) (p_1 + p_ 2 ) \text{ exp}\left ( \frac{1}{3} \text{sin} ( 8 p_1) \right )
     \bigg ]_{ k \in \{0,2,3\}} .
\end{align*}
\begin{figure}
\begin{centering}
\includegraphics[width = \textwidth]{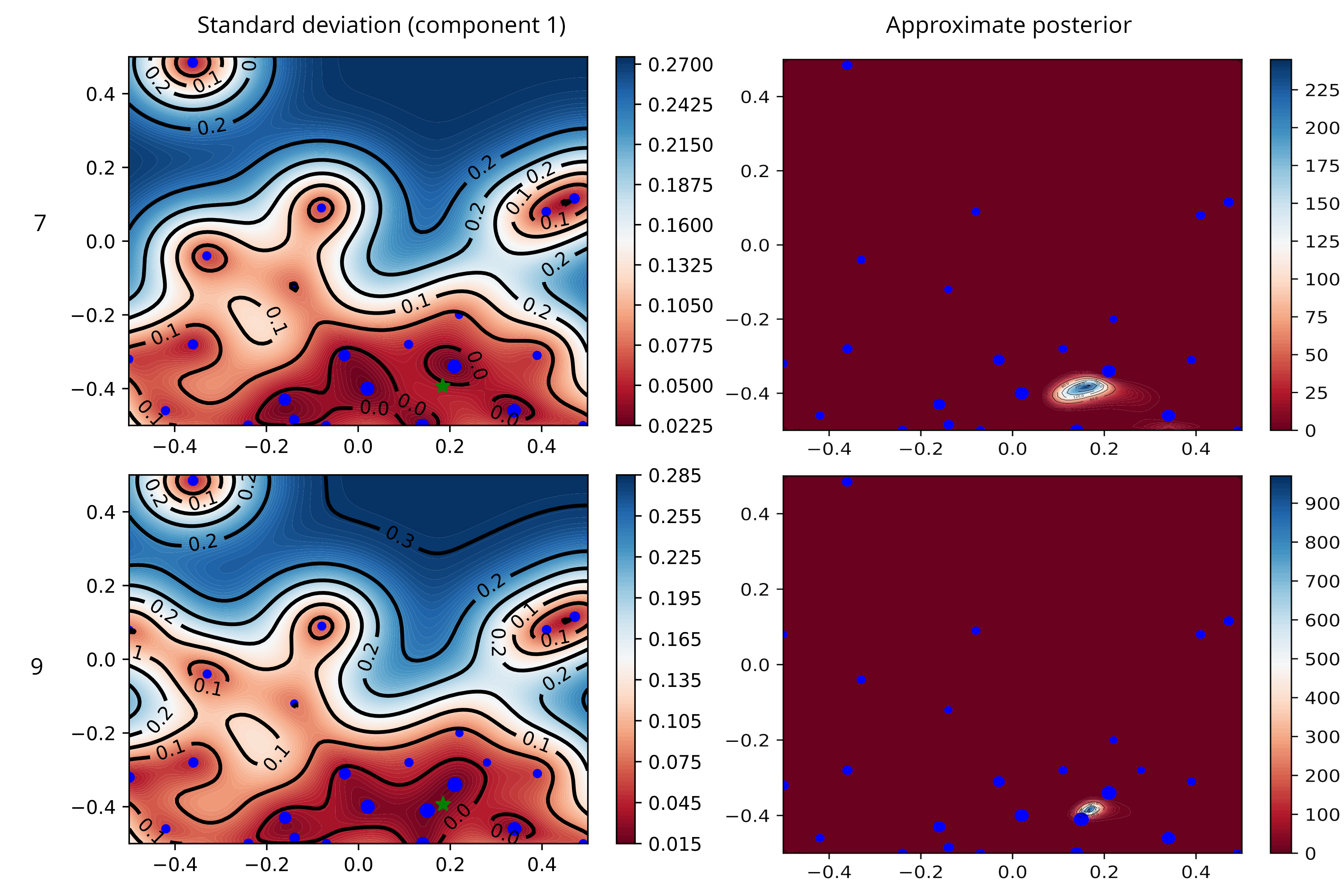}
\end{centering}
\caption{Reduction of surrogate standard deviation of the $y_1$, i.e. $k=0$, component (left) and change of posterior distribution (right) between iterations 7 and 9. The computational work for each point is represented by its size. New points are added and some of the old points are refined. The true parameter used for creating the artificial measurements $y^m$ is denoted by a green star.} 
\label{fig:standard-deviation}
\end{figure}  
The underlying model is assumed to be a quadratic FE scheme on a 3D domain, i.e. $l/r=1.5$. The discretization error is again simulated via zero mean Gaussian noise and the measurement likelihood is $\Sigma_l = 10^{-4} \mathrm{diag}( 1, 1, 4)$.

A working budget of 3600 is considered: at each iteration, 3 candidate points are considered and a fixed budget of 100 corresponding to a fixed tolerance $\tau =0.046$ is assigned to each point in the non-adaptive strategies for a total of 12 iterations. 

\begin{figure}
\begin{centering}
\includegraphics[width = \textwidth]{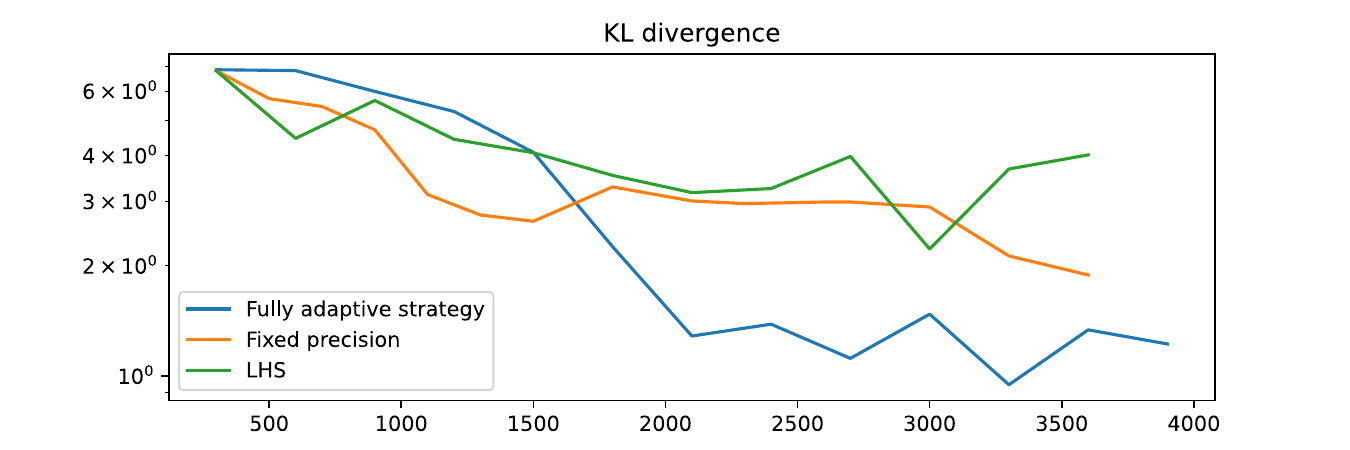}
\end{centering}
\caption{Kullback-Leibler divergence of surrogated posterior and true posterior for different training designs over the computational work spent in the 2D example.} 
\label{fig:2D}
\end{figure}

The number $n_j$ of new samples added into $\mc S_j$ is gradually increased from 200 samples at the first iteration to 4000 in the last, according to $n_j = 200 + \lfloor 26.4 j^2\rfloor$. Similarly, the number of discarded samples ranges from 200 to 2000, with $h_1 = 0$ as in the first iteration the chain is empty, and then $h_j = 200 + \lfloor 12.5 j^2\rfloor$ for $j > 1$. The error reduction by adding new points and decreasing tolerances is illustrated in Fig.~\ref{fig:standard-deviation} for a single iteration. The performance in terms of the Kullback-Leibler divergence between true and surrogated posteriors over computational work is shown in Fig.~\ref{fig:2D}. Again, a substantial performance improvement is achieved by optimizing evaluation tolerances in addition to the evaluation positions.

\section*{Conclusions}

When learning GPR surrogate models with numerically simulated training data as a replacement for the true forward model in posterior sampling, significant reductions of computational effort can be achieved with adaptive approaches. With numerical forward models that allow exploiting accuracy-work trade-offs, such as finite element simulations, the goal-oriented adaptive selection of simulation tolerances appears to be particularly effective.

\bibliographystyle{plain}
\bibliography{main}

\begin{thebibliography}{10}

\bibitem{AlvarezRosascoLawrence2012}
M.A. {\'A}lvarez, L.~Rosasco, and N.D. Lawrence.
\newblock Kernels for vector-valued functions: A review.
\newblock {\em Foundations and Trends in Machine Learning}, 4(3):195--266,
  2012.

\bibitem{BaiTeckentrupZygalakis2023}
T.~Bai, A.L. Teckentrup, and K.C. Zygalakis.
\newblock Gaussian processes for {Bayesian} inverse problems associated with
  linear partial differential equations.
\newblock Technical report, arXiv:2307.08343, 2023.

\bibitem{Crombecq}
K.~Crombecq, E.~Laermans, and T.~Dhaene.
\newblock Efficient space-filling and non-collapsing sequential design
  strategies for simulation-based modeling.
\newblock {\em European Journal of Operational Research}, 214:683--696, 2011.

\bibitem{Giunta}
A.~Giunta, S.~Wojtkiewicz, and M.~Eldred.
\newblock Overview of modern design of experiments methods for computational
  simulations (invited).
\newblock In {\em 41st Aerospace Sciences Meeting and Exhibit, AIAA 2003-649},
  pages 1--17, 2003.

\bibitem{GreenLatuszynskiPereyraRobert2015}
P.J. Green, K.~Łatuszyński, M.~Pereyra, and C.P. Robert.
\newblock Bayesian computation: a summary of the current state, and samples
  backwards and forwards.
\newblock {\em Stat. Comput.}, 25:835--862, 2015.

\bibitem{Joseph}
V.~Joseph and Y.~Hung.
\newblock Orthogonal-maximin latin hypercube designs.
\newblock {\em Statistica Sinica}, 18:171--186, 2008.

\bibitem{JarvenpaaGutmannVehtariMarttine2021}
M.~Järvenpää, M.~U. Gutmann, A.~Vehtari, and P.~Marttine.
\newblock {Parallel Gaussian process surrogate Bayesian inference with noisy
  likelihood evaluations.}
\newblock {\em Bayesian Analysis, 16, pp. 147–178.}, 2021.

\bibitem{Lehmensiek}
R.~Lehmensiek, P.~Meyer, and M.~Müller.
\newblock Adaptive sampling applied to multivariate, multiple output rational
  interpolation models with application to microwave circuits.
\newblock {\em International Journal of RF and Microwave Computer-Aided
  Engineering}, 12(4):332--340, 2002.

\bibitem{Nitzler}
J.~Nitzler, J.~Biehler, N.~Fehn, P.-S. Koutsourelakis, and A.~Wall.
\newblock A generalized probabilistic learning approach for multi-fidelity
  uncertainty quantification in complex physical simulations.
\newblock {\em Comp. Meth. Appl. Mech. Eng.}, 400:115600, 2022.

\bibitem{Queipo}
N.~Queipo, R.~Haftka, W.~Shyy, T.~Goel, R.~Vaidyanathan, and P.~Tucker.
\newblock Surrogate-based analysis and optimization.
\newblock {\em Progress in Aerospace Sciences}, 41(1):1--28, 2005.

\bibitem{RasmussenWilliams2006}
C.~Rasmussen and C.K.I. Williams.
\newblock {\em Gaussian Processes for Machine Learning}.
\newblock MIT Press, 2006.

\bibitem{SagnolHegeWeiser2016}
G.~Sagnol, H.-C. Hege, and M.~Weiser.
\newblock Using sparse kernels to design computer experiments with tunable
  precision.
\newblock In {\em Proceedings of COMPSTAT 2016}, pages 397--408, 2016.

\bibitem{SemlerWeiser2023}
P.~Semler and M.~Weiser.
\newblock Adaptive {Gaussian} process regression for efficient building of
  surrogate models in inverse problems.
\newblock {\em Inverse Problems}, 39:125003, 2023.

\bibitem{SemlerWeiser2024}
P.~Semler and M.~Weiser.
\newblock Adaptive gradient enhanced gaussian process surrogates for inverse
  problems.
\newblock In {\em Proceedings of the MATH+ Thematic Einstein Semester 2023},
  2024 (submitted).

\bibitem{SinsbeckNowak}
M.~Sinsbeck and W.~Nowak.
\newblock {Sequential Design of Computer Experiments for the Solution of
  Bayesian Inverse Problems.}
\newblock {\em SIAM/ASA Journal on Uncertainty Quantification, 5:1, 640-664.},
  2017.

\bibitem{Sugiyama}
M.~Sugiyama.
\newblock Active learning in approximately linear regression based on
  conditional expectation of generalization error.
\newblock {\em Journal of Machine Learning Research}, 7:141--–166, 2006.

\bibitem{WangBroccardo2020}
Z.~Wang and M.~Broccardo.
\newblock A novel active learning-based {Gaussian} process metamodelling
  strategy for estimating the full probability distribution in forward {UQ}
  analysis.
\newblock {\em Struct. Safety}, 84:101937, 2020.

\bibitem{WeiserGhosh2018}
M.~Weiser and S.~Ghosh.
\newblock {Theoretically optimal inexact spectral deferred correction methods}.
\newblock {\em Commu. Appl. Math. Comp. Sci.}, 13(1):53--86, 2018.

\end{thebibliography}

\end{document}